\documentclass[conference]{IEEEtran}
\IEEEoverridecommandlockouts
\usepackage{cite}
\usepackage{amsmath,amssymb,amsfonts}
\usepackage{algorithm}
\usepackage[algo2e]{algorithm2e}
\usepackage{algpseudocode}
\usepackage{graphicx}
\usepackage{textcomp}
\usepackage{xcolor}
\def\BibTeX{{\rm B\kern-.05em{\sc i\kern-.025em b}\kern-.08em
    T\kern-.1667em\lower.7ex\hbox{E}\kern-.125emX}}
\usepackage{array}
\usepackage{makecell}
\usepackage{amsfonts}
\usepackage{amssymb}
\usepackage{indentfirst}
\usepackage{amsthm}
 \usepackage{subcaption}
\usepackage{mathtools}
\usepackage{csquotes}
\usepackage{dirtytalk}
\usepackage{multirow}
\usepackage{hyperref}
\usepackage{pdflscape}
\usepackage{multirow}
\usepackage{lipsum}

\theoremstyle{definition}

\newtheorem{remark}{Remark}

\newcommand{\R}{\mathbb{R}} 
\newcommand{\J}{\mathcal{J}}

\newcommand{\prox}{\text{prox}}
\newcommand{\proj}{\text{proj}}

\newcommand{\minimize}[2]{\ensuremath{\underset{\substack{{#1}}}%
{\text{\rm minimize}}\;\;#2 }}

\begin{document}

\title{A CNC approach for Directional Total Variation\\
\thanks{This project has received funding from the European Un\-ion’s Horizon 2020 research and innovation programme under the Marie Skłodowska-Curie grant agreement No 861137.}
}



\author{\IEEEauthorblockN{Gabriele Scrivanti, \'Emilie Chouzenoux, Jean-Christophe Pesquet }
 \IEEEauthorblockA{\textit{Université Paris-Saclay, Inria, CentraleSupélec, CVN, Gif-sur-Yvette, France.} ({firstname.lastname@centralesupelec.fr)}}}
\maketitle

\begin{abstract}
The core of many approaches for the resolution of variational inverse problems arising in signal and image processing consists of promoting the sought solution to have a sparse representation in a well-suited space. A crucial task in this context is the choice of a good sparsity prior that can ensure a good trade-off between the quality of the solution and the resulting computational cost. The recently introduced Convex-Non-Convex (CNC) strategy appears as a great compromise, as it combines the high qualitative performance of non-convex sparsity-promoting functions with the convenience of dealing with convex optimization problems. This work proposes a new variational formulation to implement CNC approach in the context of image denoising. By suitably exploiting duality properties, our formulation allows to encompass sophisticated directional total variation (DTV) priors. We additionally propose an efficient optimisation strategy for the resulting convex minimisation problem. We illustrate on numerical examples the good performance of the resulting CNC-DTV method, when compared to the standard convex total variation denoiser.
\end{abstract}


\section{Introduction}
\label{sec:intro}
   A wide class of problems arising in image and signal processing can be described by sparsity-regularised variational models. The most natural sparsity inducing penalty is the $\ell_0$-pseudo-norm, but it induces the related problem to be NP hard and non-convex. A popular convex variational surrogate is the $\ell_1$-norm, though it has the drawback of under-estimating the high amplitude components of the considered signal. Non-convex variational regularisers manage to overcome this issue, at the cost of possibly introducing suboptimal local minima in the objective function. An efficient solution to keep only the best traits of these regularisers is represented by Convex-Non-Convex (CNC) strategies \cite{Selesnick2015Convex, Nikolova2015Energy,Selesnick2017Sparse,Lanza2019SParsity, Selesnick2020Nonconvex}. They consist of building convex objective functionals that include non-convex regularisation terms. Suitable CNC strategies have been designed for more and more general classes of problems in different fields of data processing (see \cite{Lanza2019SParsity} and the references therein). For the well-known Total Variation (TV) regularisation model \cite{Rudin1992Nonlinear}, a CNC modification has shown to get around the notorious problems of boundary reduction and staircasing effect for image restoration \cite{Selesnick2015Convex,Selesnick2020Nonconvex,Lanza2019SParsity} and image segmentation \cite{Chan2018ConvexNI}.
   
   In this work we propose a new variational formulation extending the CNC strategy to a popular TV-based model, namely the Directional TV (DTV). DTV was firstly introduced in \cite{Bayram2012Directional} for the restoration of images whose structures or textures follow a single dominant direction. The idea is to provide the model with a suitable directional prior information in order to allow the regularisation to operate more efficiently. The model was then extended in \cite{Zhang2013Edge} to handle multiple dominant directions by means of a spatially varying directional information that locally coincides with the edge directions of the image. In \cite{Kongskov2019DirectionalTG} and \cite{Parisotto2020Higher} the authors analysed the possibility to extend the directional approach to the Total Generalised Variation (TGV) model from  \cite{Bredies2010TotalGV}, which takes into account higher-order derivatives of the image and allows to promote piecewise-affine reconstructions, in particular. {A possibly non-convex version of DTV was proposed in \cite{calatroni2019flexible} by considering a space-variant exponent for the DTV regulariser that adapts to the local smoothness of the image.}
   
Our contributions in this work are the following. We infer a formulation for DTV that allows us to incorporate this regulariser into a CNC denoising framework. We then define a numerical procedure to efficiently address the resulting optimisation problem. We finally provide experimental results that support the use of the proposed regularisation method.

   The paper is organised as follows. In Section \ref{sec:problem_statement}, we introduce a general TV-based image denoising problem. We describe the construction of a CNC sparsity-promoting function, and provide a sufficient condition for the convexity of the overall problem. Section~\ref{sec:III} is dedicated to the derivation of our proposed CNC-DTV regularisation approach and Section~\ref{sec:IV} provides the description of a dedicated optimisation procedure to tackle it. Section~\ref{sec:results} presents our numerical results showing the interest of the proposed approach.


\section{Image denoising via CNC Total Variation}\label{sec:problem_statement}


\subsection{Notation}
Throughout this paper we denote by $\langle\cdot\mid \cdot\rangle$ the scalar product over $\R^n$ and ${\|\cdot\|_2}$ the associated Euclidean norm. $\mathsf{I}_n$ states for the $n\times n$ identity matrix. 

The {spectral norm} is defined as ${ |||L||| = \sup \{\|L z\|_{2}\mid z\in\mathbb{R}^{n}, \|z\|_{2}\leq 1 \}}$. 
 $\Gamma_0(\mathbb{R}^n)$ indicates the class of functions ${f:\mathbb{R}^n \rightarrow (-\infty,+\infty]}$ that are proper (i.e., with a nonempty domain), lower semicontinuous and convex. For a function $f\in \Gamma_0(\R^n)$, function ${f^*: \R^n \rightarrow [-\infty,+\infty] }$ represents its \textit{convex conjugate} that is defined as 
\begin{align*}
  (\forall u \in \mathbb{R}^n)\quad  f^*(u)=  \sup_{x\in \mathbb{R}^n}\{\langle x,u\rangle - f(x)\}.
\end{align*}
\subsection{Image denoising}
We focus on solving an image denoising problem, \textit{i.e.} on restoring a source image $\bar{x}\in\R^n$ from an acquired measurement $o\in \mathbb{R}^{n}$ that is related to the sought image through 
\begin{equation}
    o = \bar{x} + e,
\end{equation}
where $e \in \mathbb{R}^n$ is an additive noise, here assumed to be i.i.d. zero-mean Gaussian. Then, a simple denoising strategy consists of defining  an estimate $\hat x\in \mathbb{R}^n$ of $\bar{x}$ by solving the penalized least squares problem 
\begin{equation}
    \minimize{x\in \mathbb{R}^{n}}{ \left\{ \J(x) = \psi(x) + \frac{\lambda}{2}\|x-o\|^2_2\right\}}.
    \label{eq:denoising}
\end{equation}
Hereabove, function $\psi:\mathbb{R}^n\rightarrow (-\infty,+\infty]$ is the regularization term, associated to the regularization factor $\lambda^{-1}$, with $\lambda>0$. 


\subsection{Total-variation based regularizers}

When dealing with images, a common approach consists of choosing $\psi$ so as to sparsify the sought image in a transformed space obtained via a linear transformation. This amounts to defining $\psi  = \Psi \circ D$ where $D\in \mathbb{R}^{m\times n}$ is a linear operator and $\Psi:\mathbb{R}^m\rightarrow (-\infty,+\infty]$ is a sparsity-promoting term. Among possible choices for $\Psi$ and $D$, the total variation (TV) model \cite{Rudin1992Nonlinear} is probably the most celebrated one. TV regularizer promotes sparsity in the space of the vertical/horizontal gradients of the image, thus allowing piece-wise constant functions in the solution space. Matrix $D\in \R^{2n\times n}$ is set as the linear operator defined as $D = [D_{\rm h}^\top\ D_{\rm v}^\top]^\top$, where ${(D_{\rm h},D_{\rm v})\in (\mathbb{R}^{n\times n})^2}$ are the discrete horizontal and vertical 2D gradient operators obtained with a finite difference scheme. Then,
\begin{align}
    (\forall x\in \R^n) \quad   \operatorname{TV}(x) & = \|Dx\|_{1,2},\nonumber\\
    & = \sum_{i=1}^n \|(Dx)_i\|_2.
\end{align}
Here, for every $u\in\R^{2n}$ and  $i\in \{1,\dots,n\}$,  we use the compact notation: $\boldsymbol{u}_i = (u_{i},u_{n+i})\in \R^2$. We can also provide another definition of TV, based on duality \cite{Chambolle2004AnAF}:
\begin{align}
  \label{eq:TV:dual}
(\forall x\in&\R^n)\quad
    \operatorname{TV}(x)\nonumber\\
    &=\max_{u\in\R^{2n}} \{\langle Dx \mid u\rangle \mid \boldsymbol{u}_i \in B_2,\; i\in \{1,\dots,n\} \},
\end{align}
with $B_2 = \{ \upsilon \in \R^2 \mid  \|\upsilon \|_{2} \leq 1 \}$ the unit closed ball of $\mathbb{R}^2$.

One drawback of TV is that it gives an isotropic role to vertical and horizontal directions, that might not be well-adapted for natural images. In \cite{Bayram2012Directional}, a modified version of TV is proposed, that is more suitable for images containing objects with a (possibly non vertical/horizontal) dominant direction. The idea is to introduce an affine transformation in the dual space (i.e., the space of the image gradients). This transformation is parametrized by an expansion factor ${\alpha\geq1}$ and a rotation angle ${\theta\in[-\pi/2,\pi/2]}$. We then define the transition matrices:
\begin{equation*}
    R_{\theta} = \begin{bmatrix} \cos\theta & -\sin\theta\\ \sin\theta & \cos\theta
    \end{bmatrix},
\;\;
    \Lambda_{\alpha} = \begin{bmatrix} 1 & 0\\ 0 & \alpha
    \end{bmatrix}.
\end{equation*}
These are used to build the elliptic set $E_{\alpha,\theta}\subset \R^2$ given by $ E_{\alpha,\theta} = R_{\theta}\Lambda_{\alpha}B_2$, that is then substituted for $B_2$ in \eqref{eq:TV:dual}. 

Actually, shapes present in images might have more than one directional orientation. In order to adapt the previous regularizer to each edge directions in the image, following \cite{Zhang2013Edge}, we can rely on a set of local ellipses $E_{\alpha_i,\theta_i} \subset \R^2$ for each pixel $i\in \{1,\dots,n\}$. This leads to the so-called directional TV (DTV) regularizer, defined as
\begin{equation}
(\forall x\in\R^n)\;   \operatorname{DTV}(x;\boldsymbol\alpha,\boldsymbol\theta) = \sum_{i=1}^{n} \max_{u\in E_{{\alpha}_i,{\theta}_i}} \langle (Dx)_i| u\rangle.
\end{equation}
Here, $\boldsymbol\alpha = (\alpha_i)_{1\leq i\leq n}$, $\boldsymbol\theta= (\theta_i)_{1\leq i\leq n}$ are a predefined set of scaling factors/angles describing the edge direction information  at each pixel of the image. 

It is worth noticing that DTV is actually a generalized version of TV. Indeed, by choosing ${\alpha}_i=1$ and ${\theta}_i=0$ for every $i\in \{1,\dots,n\}$, one retrieves the classic TV functional as expressed in \eqref{eq:TV:dual}, since $E_{1,0}= R_0\Lambda_1 B_2 = B_2$.

DTV can be described in a more compact manner, by introducing the operator
${\mathcal{A}_{\boldsymbol\alpha,\boldsymbol\theta} \colon \mathbb{R}^{2n} \rightarrow \mathbb{R}^{2n}}$ that encodes the underlying local affine transformations. This linear operator is such that, for every $u \in \R^{2n}$,
  $ \mathcal{A}_{\boldsymbol\alpha,\boldsymbol\theta}u = \left(R_{\theta_i}\Lambda_{\alpha_i}\boldsymbol{u}_i\right)_{1\leq i \leq n}$. We can also express its inverse, through $ \mathcal{A}_{\boldsymbol\alpha,\boldsymbol\theta}^{-1} u = \left(\Lambda_{1/\alpha_i}R_{-\theta_i}\boldsymbol{u}_i\right)_{1\leq i \leq n}$ for every $u \in \R^{2n}$. Then, an equivalent definition of DTV is 
  \begin{align}
  \notag
  (\forall x\in&\R^n)\quad
      \operatorname{DTV}(x,\boldsymbol\alpha,\boldsymbol\theta)\\ &= \max_{u\in \mathbb{R}^{2n}}\{\langle Dx \mid u\rangle \; | \; ({A}_{\boldsymbol\alpha,\boldsymbol\theta} ^{-1}u)_i \in B_2, \; i \in \{1,\ldots,n\}\}.
      \label{eq:DTV:dual}
  \end{align}

The aforementioned (D)TV prior relies intrinsically on the convex sparsity measure $\ell_{1,2}$.
However, as emphasized in \cite{zhang2010,soubies2015continuous}, non-convexity sparsity penalties might be key to obtain high quality results. The goal of this paper is to incorporate a CNC approach with the aim to obtain an enhanced image denoiser without complexifying the optimization procedure.

\section{Proposed CNC-DTV method}
\label{sec:III}
\subsection{CNC approach}
The idea of CNC, initially proposed in \cite{Selesnick2020Nonconvex, Lanza2019SParsity}, is to use the following construction for the regularization term $\psi = \Psi_M \circ D$ in \eqref{eq:denoising}, with $D$ the discrete gradient linear operator and
\begin{equation}
    \small{
(\forall u \in \R^{2n}) \quad
    \Psi_M(u) = \varphi(u) - \underbrace{\inf_{t\in\mathbb{R}^{2n}}\left\{\varphi(t) + \frac{1}{2}\|M(u-t)\|^2_2\right\}}_{\textstyle {\varphi_M(u)}}\label{eq:ncregularizer}.}
\end{equation}
Function $\Psi_M$ is the non-convex modification of a convex sparsity promoting term $\varphi\in \Gamma_0(\R^{2n})$, parametrized by a matrix $M\in\mathbb{R}^{k\times 2n}$. It is obtained by subtracting from function $\varphi$ its so-called generalised Moreau envelope ${\varphi_M\colon\R^{2n}\rightarrow \R}$ depending on matrix $M$ {\cite[Definition 6]{Selesnick2020Nonconvex}}. This matrix plays a fundamental role, as it can be designed so as to guarantee both the efficiency of a non-convex regularisation approach and, at the same time, the overall convexity of the objective function, as we will discuss hereafter. 

In \cite{Lanza2019SParsity}, two main examples for the choice of function $\varphi$ are reported, which lead to different instances of CNC-TV regularisations, namely $\varphi =\ell_1$ yields an anisotropic version of the TV while $\varphi$ equal to the Hessian Schatten norm \cite{Lefkimmiatis2013Hessian} leads to a second-order extension of the TV regularizer. However, up to our knowledge, 
DTV (as defined in \eqref{eq:DTV:dual}) has not been explored in the context of CNC, and this is the aim of this work.

\subsection{Convexity condition}

A key feature of CNC approach is that, despite the non-convexity of the introduced penalty $\Psi_M$, it is still possible to recast the minimization problem as a convex one. More precisely, in accordance with {\cite[Example 75]{Combettes2021Fixed}}, Problem \eqref{eq:denoising} with $\psi = \Psi_M \circ D$ and $\Psi_M$ given in \eqref{eq:ncregularizer}, can be reformulated as the following minimisation problem by introducing the dual variable $y\in \mathbb{R}^{2n}$:
 \begin{align}
 \minimize{(x,y)\in \mathbb{R}^n\times \mathbb{R}^{2n} } &
    \tilde{\mathcal{J}}_M(x,y) 
     = \frac{\lambda}{2}\|x-o\|^2_2- \frac{1}{2}\|MDx\|^2_2
     \label{eq:Jm}\\
     & \notag+ \frac{1}{2}\|My\|^2_2 + \varphi(Dx)+ \varphi^{*}(M^\top M(Dx-y)).
\end{align}

\begin{remark} The equivalent problem in \eqref{eq:Jm} is slightly different from the one proposed in {\cite[Eqs (88)-(89)]{Lanza2019SParsity}} and {\cite[Eqs (52)-(53)]{Selesnick2017Sparse}}, which relied on a saddle point formulation.   
\end{remark}

{Since $\varphi$ is convex, so is its conjugate $\varphi^*$, as well as their compositions with linear mappings. It results that the only non-convex term in the objective function in \eqref{eq:Jm} is the concave quadratic one ${x\mapsto- \frac{1}{2}\|MDx\|^2_2}$}. In order to guarantee that $\tilde{\mathcal{J}}_M$
is convex, it is sufficient to ensure the 
convexity of the twice continuously differentiable function 
   ${ x\mapsto\frac{\lambda}{2}\|x-o\|^2_2- \frac{1}{2}\|MDx\|^2_2}$ \cite{Selesnick2020Nonconvex, Lanza2019SParsity,Selesnick2015Convex},
 that is equivalent to impose that the Hessian
  \begin{equation}
     H={\lambda \mathsf{I}_{2n} - D^\top M^\top M D}
     \label{eq:convcond}
 \end{equation} is a positive semidefinite matrix. In order to guarantee the existence of a unique minimizer, we further require $H$ to be positive definite and
 $M^\top M$ to be a full-rank matrix, thus making the objective function in \eqref{eq:Jm} coercive.
   

\subsection{Proposed CNC Directional Total Variation}\label{sec:proposed_DTV}

Let us now present our main contribution, that is a CNC formulation able to encompass the DTV regulariser \eqref{eq:DTV:dual} (itself being a generalisation of \eqref{eq:TV:dual}).

  First, by Fenchel-Rockafellar duality \cite{Bauschke2011ConvexAA}, from \eqref{eq:DTV:dual}, we obtain an equivalent primal formulation for DTV, by noticing that the adjoint operator
  of $\mathcal{A}^{-1}_{\boldsymbol\alpha,\boldsymbol\theta}$ is given, for every 
  $u\in \R^{2n}$, by
  $\mathcal{A}^{-*}_{\boldsymbol\alpha,\boldsymbol\theta}u=\left(R_{\theta_i}\Lambda_{1/\alpha_i}\boldsymbol{u}_i\right)_{1\leq i\leq n}$. This yields
\begin{align}
  \notag
  (\forall x\in\R^n)\quad
      &\operatorname{DTV}(x,\boldsymbol\alpha,\boldsymbol\theta)\\ &= \min_{v\in \mathbb{R}^{2n}} \{ \|v\|_{1,2} \mid \mathcal{A}^{-*}_{\boldsymbol\alpha,\boldsymbol\theta}v=Dx\}\nonumber\\
      &=  (\mathcal{A}^{-*}_{\boldsymbol\alpha,\boldsymbol\theta}\rhd \|\cdot\|_{1,2}) (Dx),
    \label{eq:DTV:primal}
\end{align}
where $(\mathcal{A}^{-*}_{\boldsymbol\alpha,\boldsymbol\theta}\rhd g)$ denotes the \textit{Exact Infimal Postcomposition} of a function $g$ by $\mathcal{A}^{-*}_{\boldsymbol\alpha,\boldsymbol\theta}$ {\cite[Chapter 12]{Bauschke2011ConvexAA}}. 
Then the proposed CNC modification of DTV parametrised by some matrix $M\in \mathbb{R}^{k\times 2n}$ reads
\begin{align*}
(\forall x\in&\mathbb{R}^n)\quad
    \operatorname{CNC-DTV}(x) \\ &=\varphi_{\mathcal{A}}(Dx) - \inf_{t\in\mathbb{R}^{2n}}\left\{\varphi_{\mathcal{A}}(t) + \frac{1}{2}\|M(Dx-t)\|^2_2\right\},\label{eq:ncdtv}
\end{align*}
where $\varphi_\mathcal{A}=\mathcal{A}^{-*}_{\boldsymbol\alpha,\boldsymbol\theta}\rhd \|\cdot\|_{1,2}$.
According to properties of the infimal postcomposition by a bounded operator {\cite[Proposition 13.24(iv)]{Bauschke2011ConvexAA}}, 
the conjugate of 
$\varphi_{\mathcal{A}}$ is expressed as 
\begin{equation}
(\forall y \in \mathbb{R}^{2n})\;\;
    \varphi^*_{\mathcal{A}} (y) = 
\iota_{B_{\infty,2}}(\mathcal{A}^{-1}_{\boldsymbol\alpha,\boldsymbol\theta}y),
\end{equation}
where $B_{\infty,2}$ is the $\ell_{\infty,2}$ ball with center 0 
and radius 1.

Hence, choosing $\varphi = \varphi_{\mathcal{A}}$ in \eqref{eq:Jm} and introducing a constraint on the range of $x$ lead to the following optimisation problem where we set ${\mathcal{H}=\mathbb{R}^n\times \mathbb{R}^{2n}\times \mathbb{R}^{2n}}$:
\begin{align}
\notag
      \minimize{(x,y,v)\in \mathcal{H} }{
      &F(x,y,v)+ \|v\|_{1,2} +\iota_{S}(x,y,v)  \\
      &+ \iota_E(Dx,\mathcal{A}^{-*}_{\boldsymbol\alpha,\boldsymbol\theta}v) 
  +\iota_{B_{\infty,2}}(\mathcal{A}^{-1}_{\boldsymbol\alpha,\boldsymbol\theta}M^{\top}M(Dx-y)).\label{eq:Jm:CNCSG:gen}}
\end{align}
Here, $S = [0,1]^n \times \R^{2n}\times \R^{2n}$ is the additional constraint set, $E$ is the vector space
$\{(w,z)\in \R^{2n}\times\R^{2n}\mid w=z\}$
and, for every $(x,y,v)\in \mathcal{H}$,
\begin{equation}
    F(x,y,v) = \frac{\lambda}{2}\|x-o\|^2_2- \frac{1}{2}\|MDx\|^2_2 + \frac{1}{2}\|My\|^2_2.
\end{equation}

{Since $\varphi_{\mathcal{A}}$ is convex, 
the positive semidefiniteness condition for
$H$ in \eqref{eq:convcond} can also be applied as a convexity guaranty for our new objective function in \eqref{eq:Jm:CNCSG:gen}.} In the next section, we investigate how to address Problem \eqref{eq:Jm:CNCSG:gen} by means of an optimisation scheme that exploits the structure of this objective function.

\section{Optimisation Algorithm}
\label{sec:IV}
\subsection{Primal-dual splitting}
The {Primal-Dual (PD) method} in \cite{condat2013Aprimaldual,Vu2013ASA,Komodakis2015Playing} allows to efficiently deal with problems involving several Lipschitzian, proximable, and linear composite terms. In \eqref{eq:Jm:CNCSG:gen}
we identify the following functions and linear operators: for every
$(y,v)\in \R^{2n}\times \R^{2n}$,
\begin{align*}
  &h_1(v) = \|v\|_{1,2} &&L_1=[ 0\; 0\;\mathsf{I}_{2n}]\\
    &h_2(v) = \iota_{B_{\infty,2}}(v) &&L_2 = [\mathcal{A}^{-1}_{\boldsymbol\alpha,\boldsymbol\theta}M^\top M D -\mathcal{A}^{-1}_{\boldsymbol\alpha,\boldsymbol\theta}M^\top M,0]\\
    &h_3(y,v) = \iota_E(y,v)&&L_3=
    \begin{bmatrix}
    D&0&0\\0&0&\mathcal{A}^{-*}_{\boldsymbol\alpha,\boldsymbol\theta}
    \end{bmatrix}, 
\end{align*}
so that \eqref{eq:Jm:CNCSG:gen} amounts to minimizing $F + \sum_{i=1}^3 h_i \circ L_i + \iota_S$.

The iterations defining a sequence ${(z_{\ell})_{\ell\in \mathbb{N}} = (x_{\ell}, y_{\ell}, v_{\ell})_{\ell\in \mathbb{N}}}$ which,  starting from point $z_0\in \mathcal{H}$, converges to a  solution ${z}_{\infty}\in \mathcal{H}$ to \eqref{eq:Jm:CNCSG:gen}, are presented in Algorithm \ref{alg:PD}, where $\delta>0$ is a Lipschitz constant of function $F$. For $f\in\Gamma_0(\R^n)$, 
$\prox_{f}$ corresponds to the {\textit{proximity operator of $f$}}.
When $f = \iota_S$, $\prox_{f}$ reduces to $\proj_{S}$, the projection onto $S$.

\begin{algorithm}[H]\caption{Primal-Dual Algorithm to solve \eqref{eq:Jm:CNCSG:gen}}
 \label{alg:PD}
 \textbf{Initialize} $z_0\in \mathcal{H}$,  $w_{1,0}\in \R^n$, $w_{2,0}\in \R^{2n}$,$w_{3,0}\in \R^{2n}$ \\
 \textbf{Set} $\tau>0$ and $\sigma>0$ s.t.
    ${\frac{1}{\tau} - \sigma|||\sum_{i=1}^3L_i^\top L_i|||\geq \frac{\delta}{2}}$\;\\
 \For{$\ell = 0,1,\ldots$}{
 {${z}_{{\ell}+1} = {\proj}_{S}(z_{\ell} - \tau\nabla F(z_{\ell}) - \tau \sum_{i=1}^3 L_i^\top w_{i,{\ell}})$}\;\\
 \For{$i=1,2,3$}{
 ${w}_{i,{\ell}+1} = \text{prox}_{\sigma h^*_i}(w_{i,{\ell}} +\sigma L_i(2{z}_{{\ell}+1} -z_{{\ell}}))$\;}
 }
\end{algorithm}

 Since $|||D|||^2 = 8$, and
 $|||\mathcal{A}^{-1}_{\boldsymbol\alpha,\boldsymbol\theta}||| \le 1$ (resp. $|||\mathcal{A}^{-*}_{\boldsymbol\alpha,\boldsymbol\theta}||| \le 1$) as $\mathcal{A}^{-1}_{\boldsymbol\alpha,\boldsymbol\theta}$ (resp. $\mathcal{A}^{-*}_{\boldsymbol\alpha,\boldsymbol\theta}$) is the combination of rotations and contractions along one axis, the norm of the involved linear operator can be upper bounded as follows:
\begin{multline}
    |||L_1^\top L_1 + L_2^\top L_2 + L_3^\top L_3||| \\
    \leq 1 + (1+8)(|||M|||^4) + 8.
\end{multline}
It is worthy to note that
Problem \eqref{eq:Jm:CNCSG:gen} and Algorithm \ref{alg:PD} represent a unified framework for convex / non-convex, classical / directional TV, since
proper choices for matrix $M$ and operator $\mathcal{A}_{\boldsymbol{\alpha},\boldsymbol{\theta}}$ allow us to model  the four different instances of the denoising problem at hand. For $M = 0_{2n} \in \R^{2n\times 2n}$, $i.e.$ the null operator, we retrieve the convex TV formulation, whereas 
 the classical TV formulation is obtained by setting, for every $i\in \{1,\ldots,n\}$,
$\alpha_i=1$ and $\theta_i= 0$.

 In the proposed framework, all the proximal computations are exact since the PD method allows us to decouple the functions that are defined by composing a convex function with a linear operator. This represents an advantage with respect to the Forward-Backward strategies proposed in {\cite[Algorithm 3]{Selesnick2020Nonconvex}} and in {\cite[Proposition 10]{Lanza2019SParsity}}, which involve nested optimisation procedures, i.e. a subroutine has to be used to compute the proximity operator of the classical TV functional.
 We dedicate the next subsection to the description of the proximity operators of the involved terms.

\subsection{Practical implementation}
The proximity operator 
of function $h_2^* = h_1 = \|\cdot\|_{1,2}$ can easily be inferred from the one of the $\ell_{2}$ norm by applying the rule for a separable sum of terms:
\begin{equation*}
   (\forall u \in \mathbb{R}^{2n}) \quad \text{prox}_{\sigma h_2^*}(u) =\left( \boldsymbol{u}_i - \frac{\boldsymbol{u}_i}{\max(\frac{\|\boldsymbol{u}_i\|_2}{\sigma},1)} \right)_{1\leq i\leq n},
\end{equation*}
while the proximity operator of 
function $h_1^*=h_2$ is obtained by applying Moreau's identity: \cite{Bauschke2011ConvexAA}
\begin{align*}
   (\forall u \in \mathbb{R}^{2n}) \quad \text{prox}_{\sigma h_1^*}(u) &= u - \sigma\prox_{\frac{h_2}{\sigma}}(\frac{u}{\sigma})\\
   &=\left(  \frac{\boldsymbol{u}_i}{\max(\|\boldsymbol{u}_i\|_2,1)} \right)_{1\leq i\leq n}.
\end{align*}
By applying again Moreau's formula, the proximity operator of $\sigma h_3^*$ is 
given by
\begin{align*}
(\forall (y,v) \in (\R^{2n})^2)\;
\text{prox} _{\sigma h_3^*}(y,v)  
    &= (y,v) - \sigma \operatorname{proj}_E
    \Big(\frac{y}{\sigma},\frac{v}{\sigma}\Big)\\
    &=\frac{1}{2} \begin{bmatrix}
    y-v\\v-y
    \end{bmatrix}.
\end{align*}
\section{Numerical Results}\label{sec:results}
We now evaluate the proposed approach for the restoration of noisy grayscale images. In accordance with \cite{Lanza2019SParsity}, we choose matrix $M = \sqrt{\gamma} \mathsf{I}_{2n}$, so that $M^\top M = \gamma\mathsf{I}_{2n}$ and therefore the convexity condition reduces to ensuring the positive semidefiniteness of ${H={\lambda\mathsf{I}_{2n}- \gamma D^\top D}}$. Since $|||D|||^2 = 8$, the convexity condition is satisfied as soon as 
$\gamma < \lambda/8$. We therefore set 
$$\gamma = (\rho\lambda)/8,$$ 
with $\rho\in\{0,0.99\}$.
For $\rho=0,$ we retrieve the convex formulation of TV, whereas for  $\rho=0.99$ we get a high degree of non-convexity for $\psi_M$ while ensuring the convexity of the global problem.

\begin{figure}
    \centering
   \resizebox{0.45\textwidth}{!}{  
    \begin{tabular}{ccc}
    \includegraphics[height=1.1in]{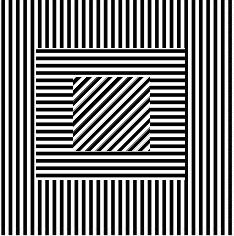}&
    \includegraphics[height=1.1in]{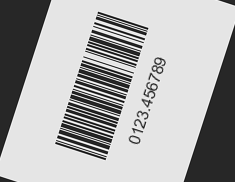}
    & \includegraphics[height=1.1in]{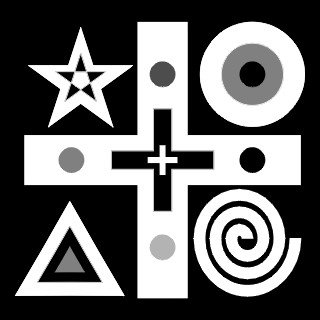}\\
    (a) & (b) & (c)\\
    \end{tabular}}
    \caption{Original images: \texttt{texture} (a), \texttt{barcode} (b), and \texttt{geometric} (c).}
    \label{fig:DTV}
\end{figure}

We illustrate the performance of CNC-DTV on three synthetic images, \texttt{texture}, \texttt{barcode}, and \texttt{geometric}, shown in Figure~ \ref{fig:DTV}. We added white zero-mean Gaussian noise with standard deviation $\sigma_e = 0.1$. In order to extract the directional information $\boldsymbol\theta$ in DTV, we exploited the strategy proposed in {\cite[Section 5]{Parisotto2020Higher}}. For the definition of parameter $\boldsymbol\alpha$, we set $\alpha_i = \alpha$ for every $i$ as in \cite{Zhang2013Edge}. We then chose the best combination of the two parameters $\alpha$ and $\lambda$ by means of a grid search to optimise the Peak Signal-to-Noise Ratio (PSNR) of the restored image.

In Table \ref{tab:dtv:psnr}, we report the best PSNR obtained when running our PD algorithm for the classic convex (i.e, $\rho = 0$) TV (C-TV), the non-convex (i.e, $\rho = 0.99$) TV (NC-TV), convex directional (C-DTV), and non-convex directional (NC-DTV) TV. This quantitative assessment shows that the combination of non-convex prior with a DTV-based space-variant regularisation yields an improvement w.r.t. the three other tested approaches. 

\begin{figure}[htb]
    \centering
    \huge{
       \resizebox{0.45\textwidth}{!}{  
    \begin{tabular}{cc}
    C-TV & NC-TV\\
    \includegraphics[trim = {0 0 0 10}, clip]{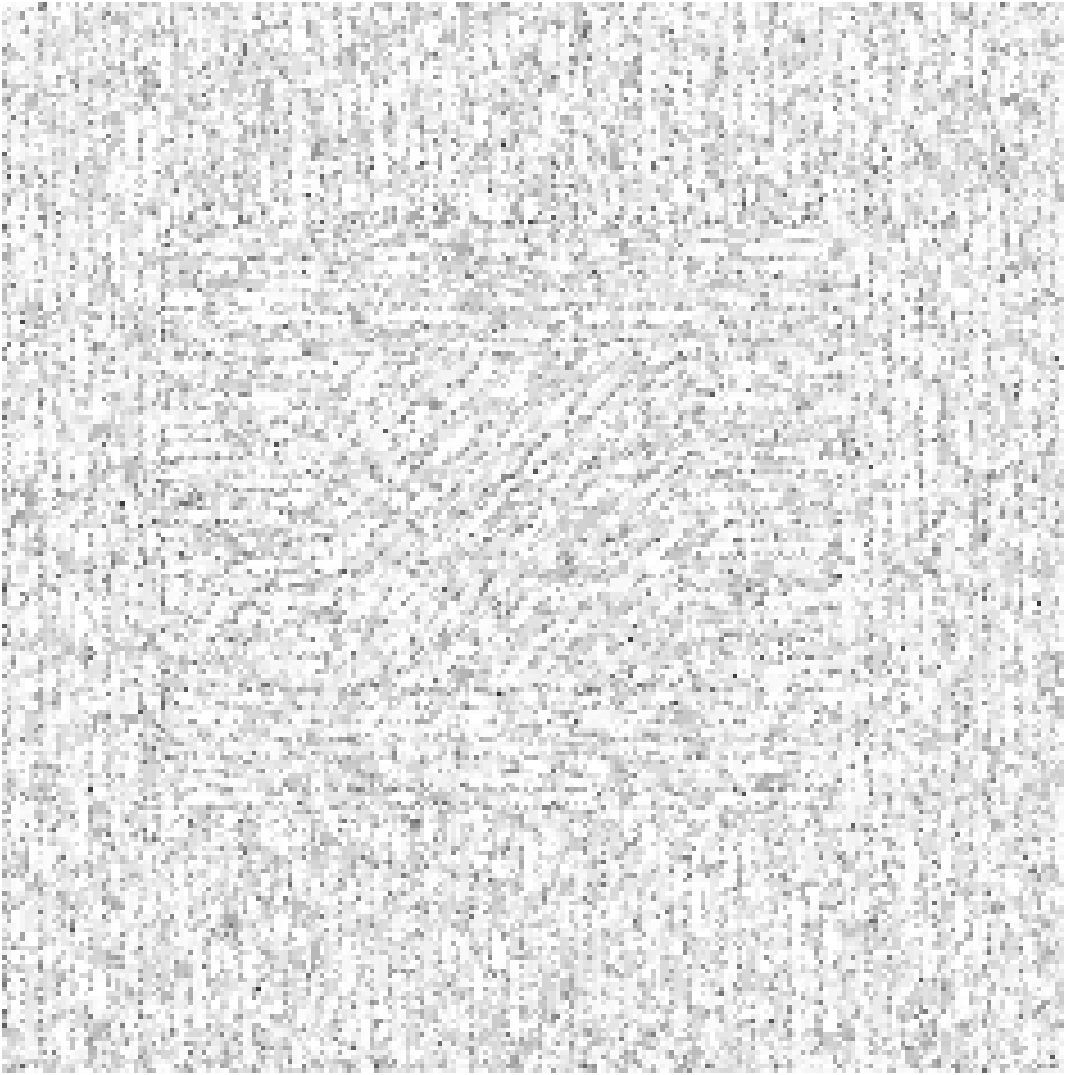}
    & \includegraphics[trim = {0 0 0 10}, clip]{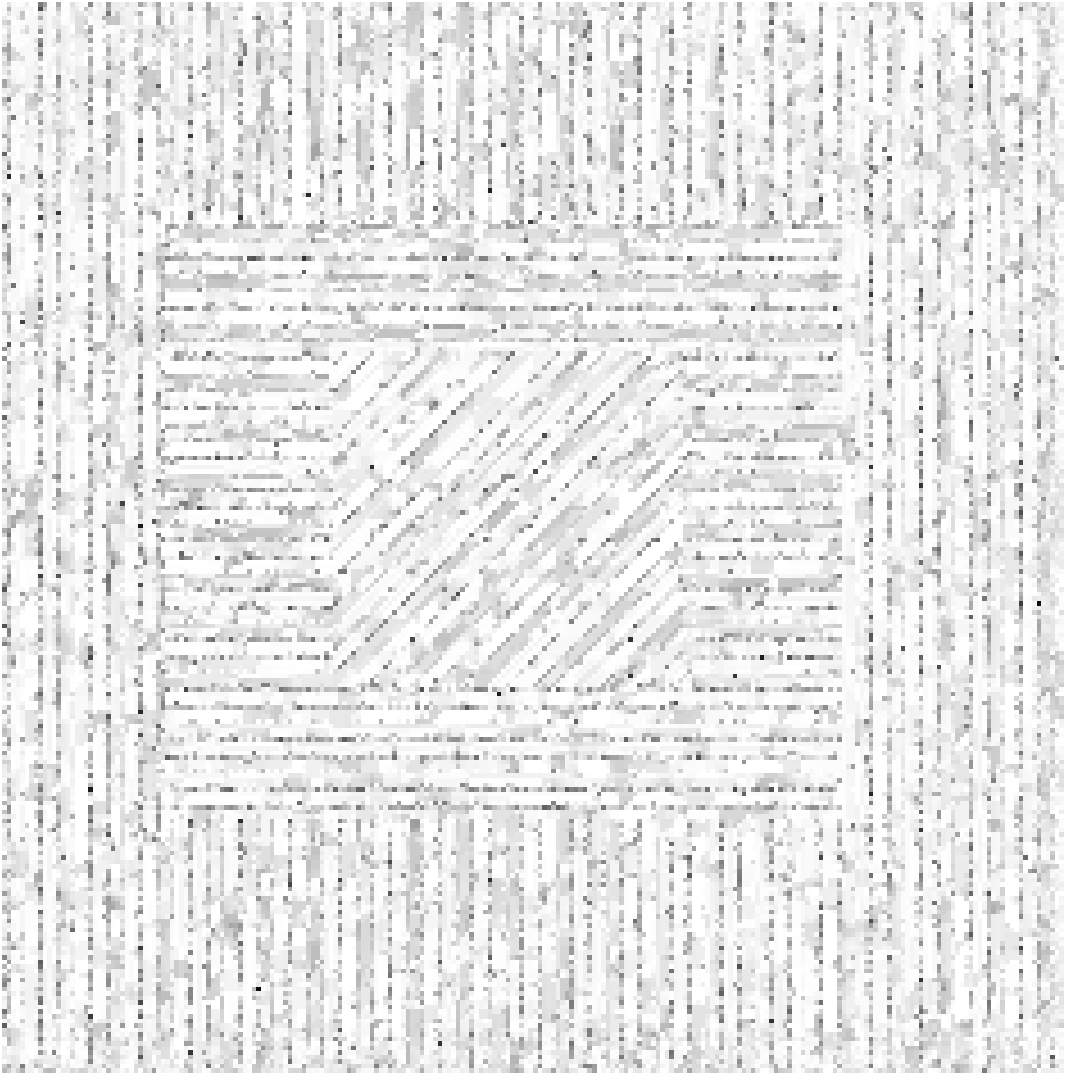}\\
    C-DTV & NC-DTV\\
    \includegraphics[trim = {0 0 0 10}, clip]{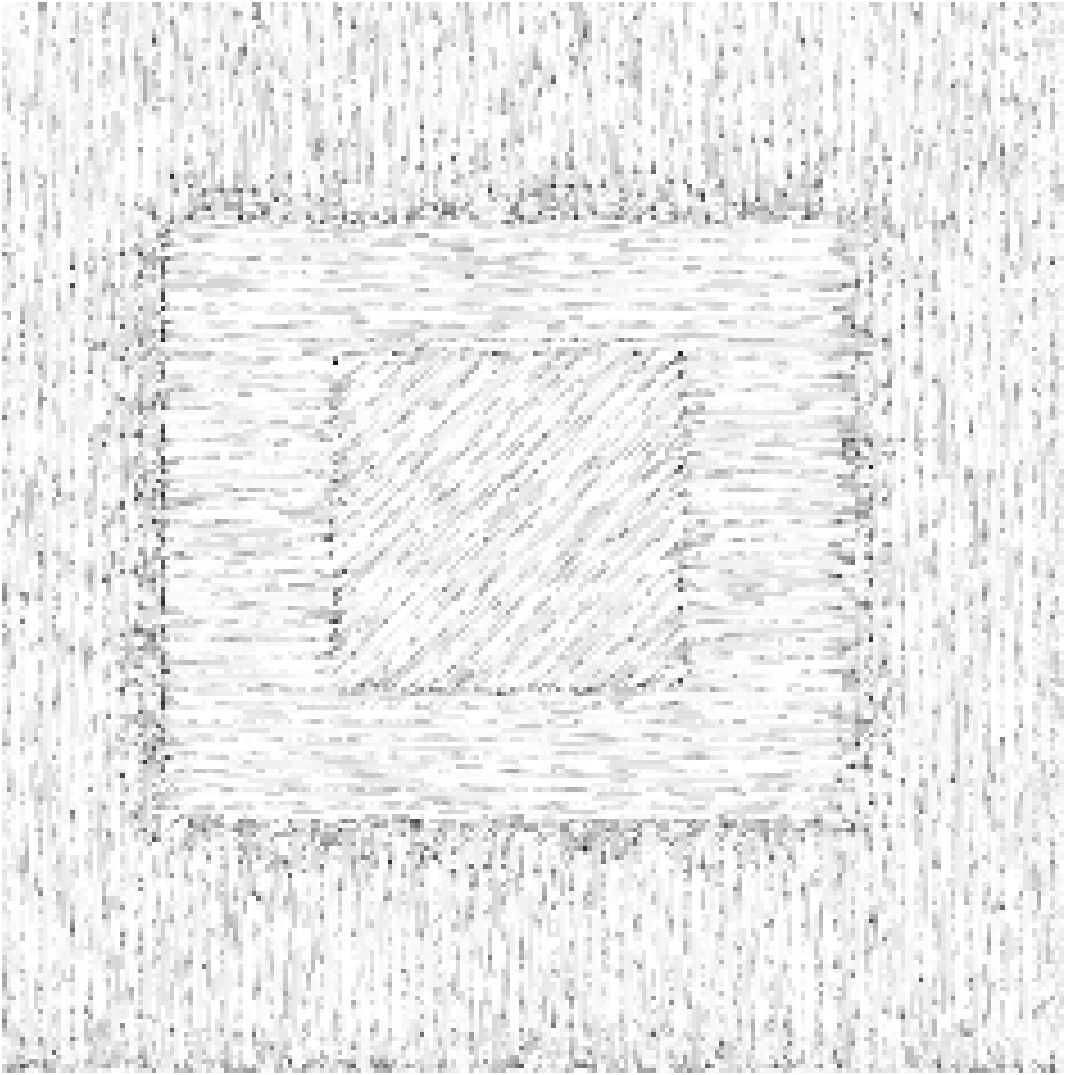}
    & \includegraphics[trim = {0 0 0 10}, clip]{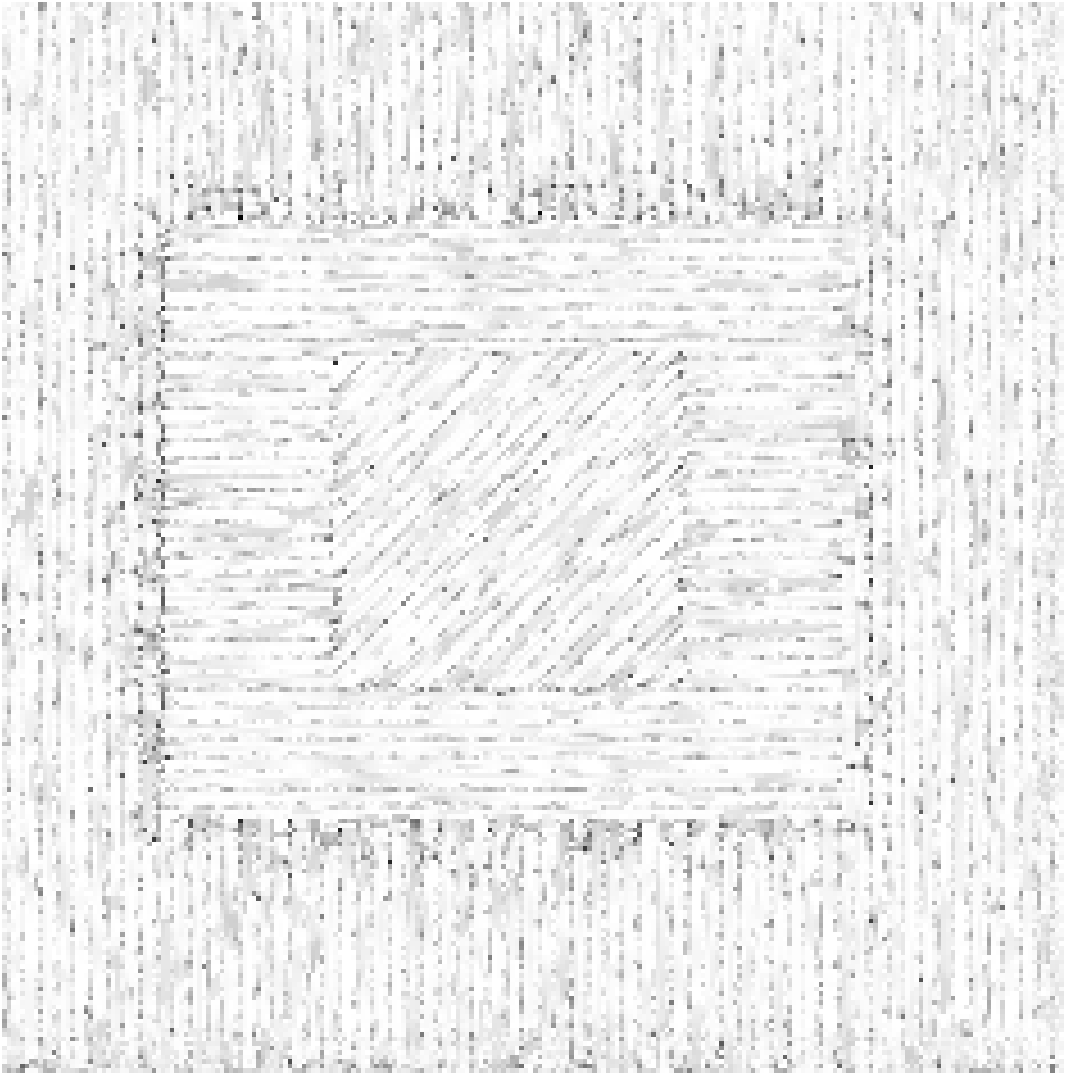}\\
    \multicolumn{2}{c}{ }\\
    \multicolumn{2}{c}{\includegraphics[width=0.5\textwidth]{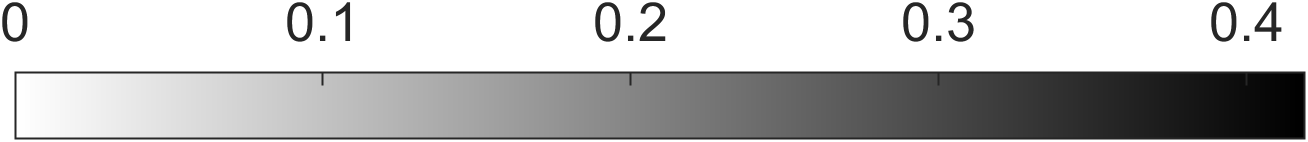}}
    \end{tabular}}}
    \caption{Residual absolute error for \texttt{texture}.}
    \label{fig:texture:absoluteerror}
\end{figure}

\begin{table}[h]
    \centering
    \resizebox{0.45\textwidth}{!}{  
    \begin{tabular}{|c||cccc|}
    \hline
       image     &   C-TV  &  CNC-TV &   C-DTV ($1/\alpha$)  &  CNC-DTV ($1/\alpha$)\\
            \hline\hline
   \texttt{texture}         &  23.70   &  24.34    & 25.83 (0.1)   &  \textbf{26.76} (0.25) \\
        \texttt{barcode}  &  25.60   &  26.92    & 26.07 (0.4)   &  \textbf{27.69} (0.5) \\
    \texttt{geometric}    &  30.81   &   31.16   &   31.24 (0.45)   &  \textbf{32.30} (0.45)\\
       \hline
    \end{tabular}}
    \caption{Best PSNR for $\sigma_e=0.1$ optimised over a grid search for parameters $\lambda$ and $\alpha$.}
    \label{tab:dtv:psnr}
\end{table}
We also provide in Figure \ref{fig:texture:absoluteerror} a visual illustration of the performance of the four approaches on \texttt{texture} image, showing the absolute error between the best estimated solution $\hat{x}$ and the original image $\bar{x}$.
The reconstructions involving a non-convex regulariser or directional information show peculiar structures in the distribution of its absolute residual error, whereas the classic convex approach yields a rather dense error distribution. NC-DTV inherits the high coherence w.r.t the directions in the image from C-DTV and the high accuracy in noise removal and sharp transition reconstruction from NC-TV.

\begin{figure}[htb]
    \centering
    \resizebox{0.49\textwidth}{!}{  
    \begin{tabular}{cccc}
      \makecell{\rotatebox{90}{\huge{PSNR} }} & \makecell{\includegraphics{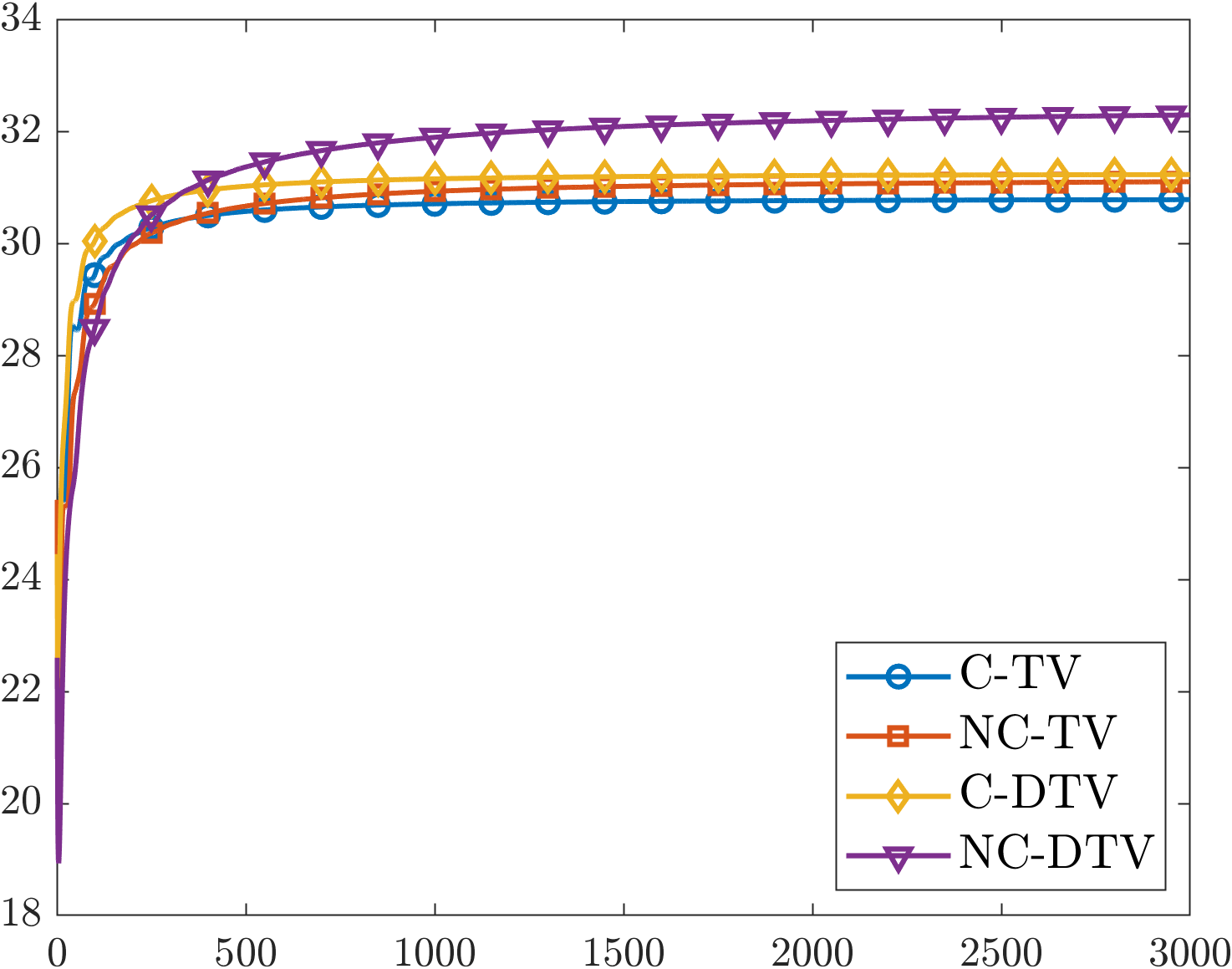}} &\makecell{\rotatebox{90}{\huge{$\|z_{\ell} - z_{\infty}\|_2/\|z_{\infty}\|_2$} }} &
      \makecell{\includegraphics{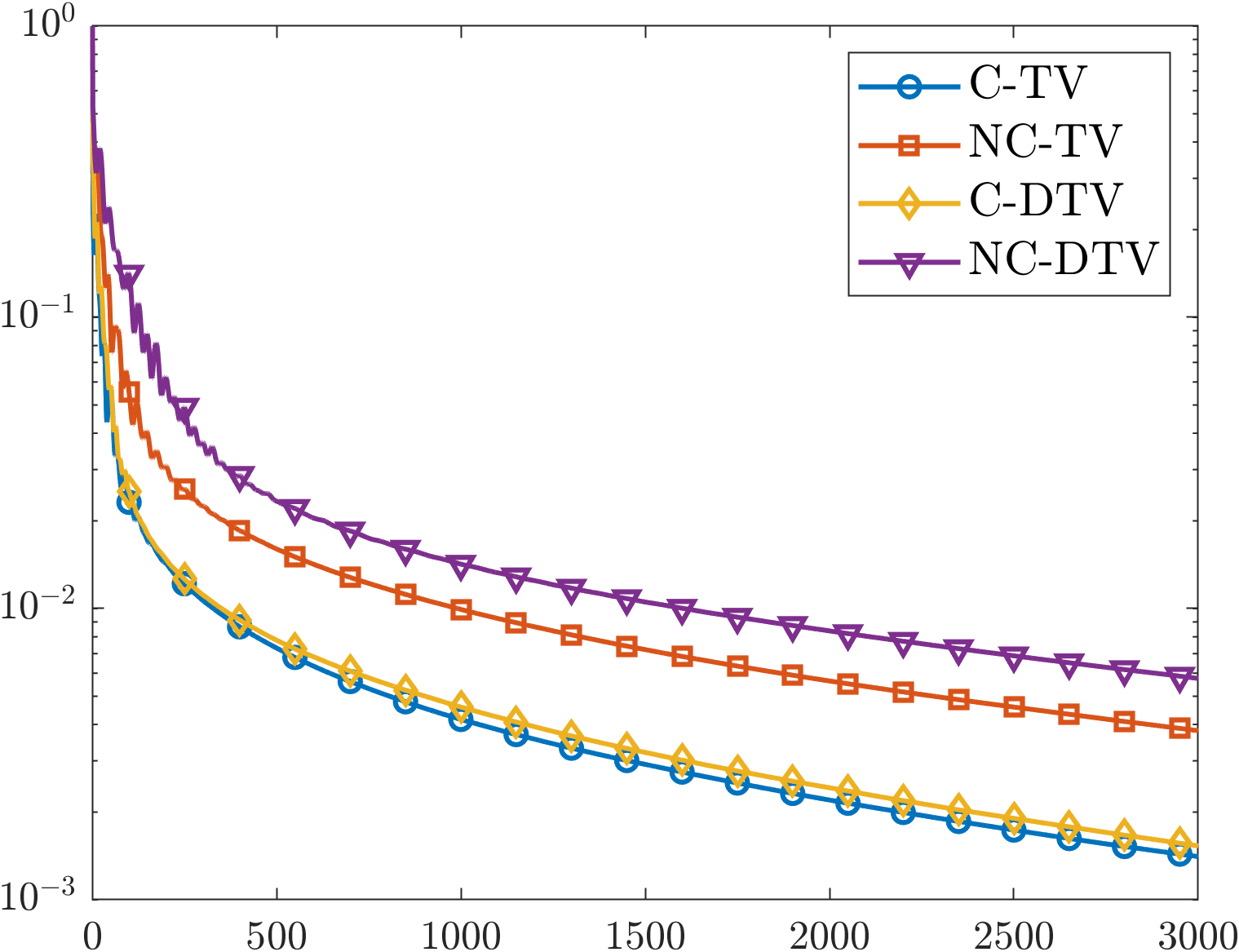} }\\
      & \huge{$\ell$} & &\huge{$\ell$}\\
    \end{tabular}}
    \caption{PSNR versus iterations (left) and distance from the iterates
to the solution versus iterations in logarithmic scale (right) for \texttt{geometric}.}
    \label{fig:convergence}
\end{figure}

Eventually, in Figure \ref{fig:convergence}, we illustrate for image \texttt{geometric} the PSNR evolution along 3000 iterations (left) and the distance from the iterates $z_\ell$ to the solution $z_\infty$ in logarithmic scale (right) of the four approaches, which shows the fast and stable convergence behaviour of the proposed PD algorithm.
      
\section{Conclusions}
In this work we investigated the extension of the CNC approach to a directional version of the TV regularisation model for image denoising. We proposed to address the resulting minimisation problem with a primal-dual procedure that efficiently exploits the structure of the objective function and we provided numerical results supporting the interest of the proposed approach.
\bibliographystyle{IEEEtran}
\bibliography{strings,refs}

\end{document}